\documentclass[11pt]{article}

\usepackage{amssymb,amsmath,amscd,amsthm,dsfont,xy,enumitem,mathrsfs,pstricks}
\xyoption{all}
\theoremstyle{definition}

\theoremstyle{remark}

\def\lra{\longrightarrow}

\def\BE#1{\begin{equation}\label{#1}}
\def\EE{\end{equation}}
\def\eref#1{(\ref{#1})}

\def\tn#1{\textnormal{#1}}

\def\ov#1{\overline{#1}}
\def\tn#1{\textnormal{#1}}

\def\sm#1{\begin{small}#1\end{small}}

\def\C{\mathbb C}
\def\cC{\mathcal C}

\def\cM{\mathcal M}

\def\cO{\mathcal O}
\def\P{\mathbb P}

\def\Q{\mathbb Q}
\def\R{\mathbb R}

\def\Z{\mathbb Z}

\def\io{\iota}

\def\om{\omega}
\def\si{\sigma}

\def\vr{\varrho}

\def\Si{\Sigma}

\def\eset{\emptyset}
\def\i{\infty}

\def\1{\mathds{1}}

\begin{document}

\title{The (Co)homology of the Deligne-Mumford Moduli Spaces\\ of Marked Rational Curves} 
\author{Aleksey Zinger}
\date{\today}
\maketitle

\begin{abstract}
\noindent
This informal note collects key results and open problems on the (co)homology 
of the Deligne-Mumford moduli spaces of real marked rational curves.
The open problems are both of topological nature, aiming to investigate the (co)homology 
of these spaces further, and of algebraic nature, aiming to describe the structure
of the homology of these spaces in operad-like terms.\\
\\ 
\end{abstract}

\noindent
$\ov\cM_{\ell}=$ Deligne-Mumford moduli space of (complex) rational curves 
with $\ell$~marked points; smooth projective variety of complex dimension~$\ell\!-\!3$.
An element of this space is the equivalence class of stable tuples
$$\cC\equiv\big(\Si,z_1,\ldots,z_{\ell}\big),$$
where $\Si$ is a connected, possibly nodal, Riemann surface and $z_1,\ldots,z_{\ell}\!\in\!\Si$
are distinct smooth points.\\

\noindent
$\R\ov\cM_{k,\ell}=$ Deligne-Mumford moduli space of real rational curves 
with $k$~real marked points and $\ell$ conjugate pairs of marked points;
compact smooth real analytic variety of dimension~$k\!+\!2\ell\!-\!3$;
orientable if $k\!=\!0$ or $k\!+\!2\ell\!\le\!4$.
An element of this space is the equivalence class of stable tuples
$$\cC\equiv\big(\Si,\si,z_1,\ldots,z_k,(z_1^+,z_1^-),\ldots,(z_{\ell}^+,z_{\ell}^-)\!\big),$$
where $\Si$ is a connected, possibly nodal, Riemann surface, 
$\si$ is an anti-holomorphic involution on~$\Si$ and 
$$z_1,\ldots,z_k,z_1^+,z_1^-,\ldots,z_{\ell}^+,z_{\ell}^-\in\Si$$
are distinct smooth points such that $\si(z_i^+)\!=\!z_i^-$ for $i\!=\!1,\ldots,\ell$.\\

\noindent
$[\ell]\!=\!\{1,\ldots,\ell\}$. 

\section{Key Results}

\begin{figure}
\begin{pspicture}(-2.8,-1.6)(10,1.2)
\psset{unit=.3cm}
% left diagram
\pscircle*(5,0){.2}
\psline[linewidth=.07](4,-1)(8,3)\psline[linewidth=.07](4,1)(8,-3)
\pscircle*(6,1){.2}\pscircle*(7,2){.2}\pscircle*(6,-1){.2}\pscircle*(7,-2){.2}
\rput(6,2){\sm{$1$}}\rput(7,3){\sm{$2$}}\rput(6,-2){\sm{$3$}}\rput(7,-3){\sm{$4$}}
\rput(12,0){\begin{Large}$=$\end{Large}}
\rput(6,-5){$D_{\{1,2\},\{3,4\}}$}
% middle diagram
\pscircle*(17,0){.2}
\psline[linewidth=.07](16,-1)(20,3)\psline[linewidth=.07](16,1)(20,-3)
\pscircle*(18,1){.2}\pscircle*(19,2){.2}\pscircle*(18,-1){.2}\pscircle*(19,-2){.2}
\rput(18,2){\sm{$1$}}\rput(19,3){\sm{$3$}}\rput(18,-2){\sm{$2$}}\rput(19,-3){\sm{$4$}}
\rput(24,0){\begin{Large}$=$\end{Large}}
\rput(18,-5){$D_{\{1,3\},\{2,4\}}$}
% middle diagram
\pscircle*(29,0){.2}
\psline[linewidth=.07](28,-1)(32,3)\psline[linewidth=.07](28,1)(32,-3)
\pscircle*(30,1){.2}\pscircle*(31,2){.2}\pscircle*(30,-1){.2}\pscircle*(31,-2){.2}
\rput(30,2){\sm{$1$}}\rput(31,3){\sm{$4$}}\rput(30,-2){\sm{$2$}}\rput(31,-3){\sm{$3$}}
\rput(30,-5){$D_{\{1,4\},\{2,3\}}$}
\end{pspicture}
\caption{Relations between the divisors $D_{J,K}\!\subset\!\ov\cM_4$ in $H^2(\ov\cM_4;\Z)$.}
\label{M04rel_fig}
\end{figure}
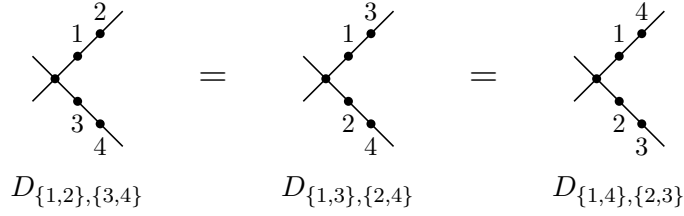

\noindent
{\bf Complex case, $\ov\cM_{\ell}$.} This smooth projective variety contains smooth divisors~$D_{J,K}$,
with $J\!\sqcup\!K\!=\![\ell]$ and $|J|,|K|\!\ge\!2$, 
whose generic element is a two-component curve with the marked points
on the two components indexed by~$J$ and~$K$; see Figure~\ref{M04rel_fig} for some examples. 
It is a now classical result of Keel~\cite{Keel} that 
\begin{enumerate}[label=($\C\arabic*$),ref=\arabic*,leftmargin=*]

\item\label{Cgen_it} these divisors generate  $H^*(\ov\cM_{\ell};\Z)$ as an algebra,

\item\label{Crel_it} subject only to some obvious relations; 

\end{enumerate}
see \cite[Theorem~2.1]{RDMhomol} for a precise statement.
The relations arise from
\begin{enumerate}[label=($\C\ref{Crel_it}\alph*$),leftmargin=*]

\item certain divisors~$D_{J,K}\!\subset\!\ov\cM_{\ell}$ being disjoint and

\item all divisors~$D_{J,K}$ in~$\ov\cM_4\!\approx\!\C\P^1$ being points 
and thus equivalent in $H^2(\ov\cM_4;\Z)$, as indicated in Figure~\ref{M04rel_fig}.

\end{enumerate}
The above relations in $H^2(\ov\cM_4;\Z)$ pull back to relations in  $H^2(\ov\cM_{\ell};\Z)$
with $\ell\!\ge\!4$ via the forgetful morphisms $\ov\cM_{\ell}\!\lra\!\ov\cM_4$ dropping all but
four of the marked points (every 4-element subset of~$[\ell]$ determines such a morphism).\\

\noindent
The proof in~\cite{Keel} presents $\ov\cM_{\ell+1}$ as a sequence of holomorphic blowups 
of~$\ov\cM_{\ell}\!\times\!\ov\cM_4$ and determines $H^*(X_{\vr};\Z)$ inductively for 
every intermediate blowup~$X_{\vr}$, assuming the claimed result for~$H^*(\ov\cM_{\ell};\Z)$.
The proof of the main result of~\cite{Keel} given in~\cite{RDMhomol} still uses this sequence
of blowups to obtain~(CDM\ref{Cgen_it}), but establishes~(CDM\ref{Crel_it}) directly,
without determining $H^*(X_{\vr};\Z)$ for the intermediate blowups~$X_{\vr}$.\\

\noindent
For any commutative ring $R$ with unity~1, the modules $H_*(\ov\cM_{\ell};R)$ can be arranged into
an operad of graded $R$-modules
$$\cO_{\C}(\ell)\equiv\begin{cases}R,&\hbox{if}~\ell\!=\!1;\\
H_*(\ov\cM_{\ell+1};R),&\hbox{if}~\ell\!\ge\!2;\end{cases}$$
as follows.
For integers $k,\ell\!\ge\!2$ and $i\!\in\![k]$, let
$$\io_{k\ell i}\!: \ov\cM_{k+1}\!\times\!\ov\cM_{\ell+1}\lra\ov\cM_{k+\ell}$$
be the immersion identifying the $i$-th marked point of a marked curve 
$[\cC_1]\!\in\!\ov\cM_{k+1}$ with the last marked point of a marked curve 
$[\cC_2]\!\in\!\ov\cM_{\ell+1}$, keeping the indices of the first $i\!-\!1$ marked points of~$[\cC_1]$
the same,
increasing the indices of the last $k\!+\!1\!-\!i$ marked points of~$[\cC_1]$
by~$\ell\!-\!1$,
and increasing the indices of the first $\ell$~marked points of~$[\cC_2]$ by~$i\!-\!1$.
This induces homomorphisms
$$\io_{k\ell i*}\!: \cO_{\C}(k)\!\times\!\cO_{\C}(\ell)\!\approx\!
H_*\big(\ov\cM_{k+1}\!\times\!\ov\cM_{\ell+1};R\big)\lra\cO_{\C}(k\!+\!\ell\!-\!1).$$
For $k,\ell\!\in\!\Z^+$ and $i\!\in\![k]$, define
\BE{prodcases_e}\begin{aligned}
\io_{k1i*}\!: \cO_{\C}(k)\!\times\!\cO_{\C}(1)&\lra\cO_{\C}(k), &\qquad
\io_{k1i*}(\eta,r)&=r\eta,\\
\io_{1\ell 1*}\!: \cO_{\C}(1)\!\times\!\cO_{\C}(\ell)&\lra\cO_{\C}(\ell), &\qquad
\io_{1\ell 1*}(r,\eta)&=r\eta.
\end{aligned}\EE
For $k,\ell_1,\ldots,\ell_k\!\in\!\Z^+$, define
\begin{gather}\label{expandcomp_e}
\circ\!:\cO_{\C}(k)\!\times\!\cO_{\C}(\ell_1)\!\times\!\ldots\!\times\!\cO_{\C}(\ell_k)\lra
\cO_{\C}\big(\ell_1\!+\!\ldots\!+\!\ell_k\big),\\
\notag
\eta\!\circ\!(\eta_1,\ldots,\eta_k)
=\io_{(1+\ell_2+\ldots+\ell_k)\ell_11*}
\big(\ldots\io_{(k-1+\ell_k)\ell_{k-1}(k-1)*}\big(\io_{k\ell_k k*}(\eta,\eta_k),\eta_{k-1}\big),
\ldots,\eta_1\big). 
\end{gather}
The resulting operad $(\!(\cO_{\C}(\ell))_{\ell\in\Z^+},\circ)$ is
 the operad of hypercommutative algebras; see~\cite{Getz95}.\\

\noindent
{\bf Real case with real points only, $\R\ov\cM_{k,0}$.}
By the main result of~\cite{Krasnov}, $H^*(\R\ov\cM_{k,0};\Z_2)$ is described
in the same way as $H^{2*}(\ov\cM_k;\Z)$ above with the complex subvarieties 
$D_{J,K}\!\subset\!\ov\cM_k$ replaced by the analogous real subvarieties 
\hbox{$\R D_{J,K}\!\subset\!\R\ov\cM_{0,k}$}; see the proof of \cite[Theorem~5.6]{EHKR}.
Since $\R\ov\cM_{k,0}$ is not orientable for $k\!\ge\!5$,
the homology and cohomology of~$\R\ov\cM_{k,0}$ with coefficients in a ring of characteristic
larger than~2 are quite distinct.
For any commutative ring~$R$ with unity and distinct \hbox{$a,b,c,d\!\in\![k]$}, 
$H^1(\R\ov\cM_{k,0};R)$ contains the pullback~$\om_{abcd}$ of the Poincare dual
of a point in $\R\ov\cM_{4,0}\!\approx\!\R\P^1$ by the forgetful morphism
\hbox{$\R\ov\cM_{k,0}\!\lra\!\R\ov\cM_{4,0}$} sending the marked indexed by $a,b,c,d$
to the marked points indexed by~$1,2,3,4$, respectively, and dropping the remaining marked points.
This class is Poincare dual to a linear combination of the real hypersurfaces~$\R D_{J,K}$.\\

\noindent
By \cite[Theorem~2.9]{EHKR},
\begin{enumerate}[label=($\R\R\arabic*$),ref=\arabic*,leftmargin=*]

\item\label{RRgen_it} the classes~$\om_{abcd}$ generate  $H^*(\R\ov\cM_{k,0};\Q)$ as an algebra,

\item\label{RRrel_it} subject only to some expected relations.

\end{enumerate}
The relations arise from
\begin{enumerate}[label=($\R\R\ref{RRrel_it}\alph*)$,leftmargin=*]

\stepcounter{enumi}

\item the diffeomorphism of $\R\ov\cM_{4,0}$ induced by the interchange of two marked points
being orientation-reversing and 

\item $\R\ov\cM_{5,0}$ being diffeomorphic to the blowup of $\R\P^1\!\times\!\R\P^1$
at three points of the diagonal.

\end{enumerate}
The resulting relations in $H^1(\R\ov\cM_{4,0};\Q)$, $H^1(\R\ov\cM_{5,0};\Q)$, and $H^2(\R\ov\cM_{5,0};\Q)$
pull back to relations in $H^1(\R\ov\cM_{k,0};\Q)$ and $H^2(\R\ov\cM_{k,0};\Q)$ via 
the forgetful morphisms 
$$\R\ov\cM_{k,0}\lra \R\ov\cM_{4,0}  \qquad\hbox{and}\qquad
\R\ov\cM_{k,0}\lra \R\ov\cM_{5,0}$$
dropping all but four (resp.~five) of the marked points.
By \cite[Corollary~6.9]{EHKR}, certain natural subvarieties of~$\R\ov\cM_{k,0}$
form a linear basis for~$H_*(\R\ov\cM_{k,0};\Q)$.
By \cite[Corollary~3.8]{Rains}, $H_*(\R\ov\cM_{k,0};\Z)$ and $H^*(\R\ov\cM_{k,0};\Z)$
have only 2-torsion; it is described in \cite[Section~5]{EHKR}.\\

\noindent
Similarly to the complex case, there are node-identifying immersions
$$\io_{k\ell i}^{\R}\!: \R\ov\cM_{k+1,0}\!\times\!\R\ov\cM_{\ell+1,0}\lra\R\ov\cM_{k+\ell,0}.$$
They again determine an operad structure on the sequence
$$\cO_{\R\R}(k)\equiv\begin{cases}\Q,&\hbox{if}~k\!=\!1;\\
H_*(\R\ov\cM_{k+1,0};\Q),&\hbox{if}~k\!\ge\!2;\end{cases}$$
of graded vector spaces.
By \cite[Theorem~2.14]{EHKR},
this operad is the operad of unital 2-Gerstenhaber algebras.\\

\noindent
The proof of \cite[Corollary~3.8]{Rains} is based on presenting~$\R\ov\cM_{k,0}$ 
as a sequence of real blowups of~$(\R\P^1)^{k-3}$ and shows that 
$H_*(X_{\vr};\Z)$ has only 2-torsion for every intermediate blowup~$X_{\vr}$.
Theorems~2.9 and~2.12 and Corollary~6.9 in~\cite{EHKR} are proved in parallel,
by combining the aforementioned topological results 
of \cite{Keel} and~\cite{Krasnov} with algebraic considerations.
The same results should be obtainable by presenting~$\R\ov\cM_{k+1,0}$ 
as a sequence of real blowups of~$\R\ov\cM_{k,0}\!\times\!\R\ov\cM_{4,0}$ 
as in~\cite{RDMbl} and relating the (co)homologies of the intermediate blowups
as in~\cite[Section~3.2]{RDMhomol}.\\

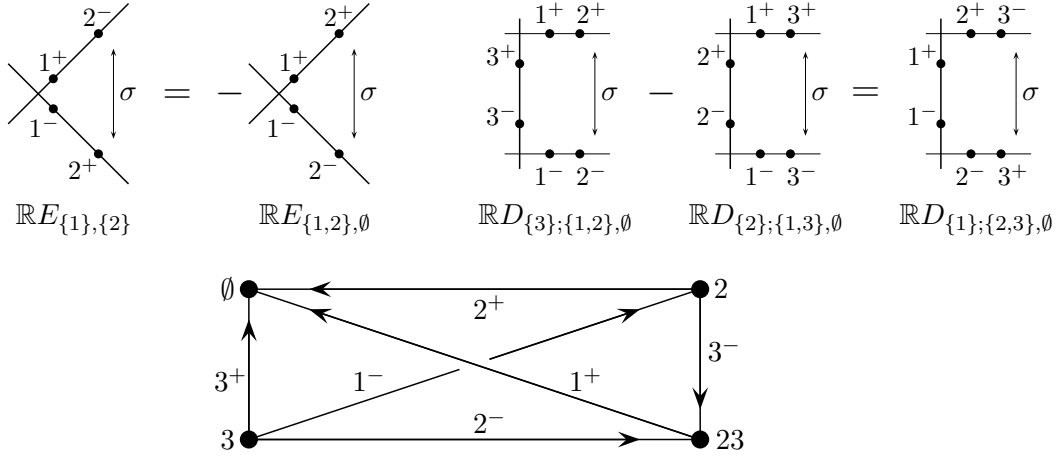
\begin{figure}
\begin{pspicture}(-.5,-1.4)(10,5.4)
\psset{unit=.4cm}
% M02 eqn, 1st diagram
\psline[linewidth=.05](2,9)(6,13)\psline[linewidth=.05](2,11)(6,7)
\pscircle*(3.5,10.5){.15}\pscircle*(5,12){.15}
\pscircle*(3.5,9.5){.15}\pscircle*(5,8){.15}
\rput(3.2,8.9){\sm{$1^-$}}\rput(4.5,7.5){\sm{$2^+$}}
\rput(3.5,11.1){\sm{$1^+$}}\rput(5,12.6){\sm{$2^-$}}
\psline[linewidth=.025]{<->}(5.5,8.5)(5.5,11.5)\rput(6,10){$\si$}
\rput(8.5,10){\begin{Large}$=~-$\end{Large}}
\rput(4.2,5.8){$\R E_{\{1\},\{2\}}$}
% M02 eqn, 2nd diagram
\psline[linewidth=.05](10,9)(14,13)\psline[linewidth=.05](10,11)(14,7)
\pscircle*(11.5,10.5){.15}\pscircle*(13,12){.15}
\pscircle*(11.5,9.5){.15}\pscircle*(13,8){.15}
\rput(11.1,8.9){\sm{$1^-$}}\rput(12.5,7.5){\sm{$2^-$}}
\rput(11.5,11.1){\sm{$1^+$}}\rput(13,12.6){\sm{$2^+$}}
\psline[linewidth=.025]{<->}(13.5,8.5)(13.5,11.5)\rput(14,10){$\si$}
\rput(12.2,5.8){$\R E_{\{1,2\},\eset}$}
% M03 eqn, 1st diagram
\psline[linewidth=.05](19,12.5)(19,7.5)
\psline[linewidth=.02](18.5,12)(22,12)\psline[linewidth=.02](18.5,8)(22,8)
\pscircle*(20,12){.15}\pscircle*(21,12){.15}
\pscircle*(19,11){.15}\pscircle*(19,9){.15}
\rput(20,12.7){\sm{$1^+$}}\rput(21.4,12.7){\sm{$2^+$}}
\rput(18.4,11.3){\sm{$3^+$}}\rput(18.4,9.3){\sm{$3^-$}}
\pscircle*(20,8){.15}\pscircle*(21,8){.15}
\rput(20,7.3){\sm{$1^-$}}\rput(21.4,7.3){\sm{$2^-$}}
\psline[linewidth=.025]{<->}(21.5,8.6)(21.5,11.4)\rput(22,10){$\si$}
\rput(23.7,10){\begin{Large}$-$\end{Large}}
\rput(20.2,5.8){$\R D_{\{3\};\{1,2\},\eset}$}
% M03 eqn, 2nd diagram
\psline[linewidth=.05](26,12.5)(26,7.5)
\psline[linewidth=.02](25.5,12)(29,12)\psline[linewidth=.02](25.5,8)(29,8)
\pscircle*(27,12){.15}\pscircle*(28,12){.15}
\pscircle*(26,11){.15}\pscircle*(26,9){.15}
\rput(27,12.7){\sm{$1^+$}}\rput(28.4,12.7){\sm{$3^+$}}
\rput(25.4,11.3){\sm{$2^+$}}\rput(25.4,9.3){\sm{$2^-$}}
\pscircle*(27,8){.15}\pscircle*(28,8){.15}
\rput(27,7.3){\sm{$1^-$}}\rput(28.4,7.3){\sm{$3^-$}}
\psline[linewidth=.025]{<->}(28.5,8.6)(28.5,11.4)\rput(29,10){$\si$}
\rput(30.5,10){\begin{Large}$=$\end{Large}}
\rput(27.2,5.8){$\R D_{\{2\};\{1,3\},\eset}$}
% M03 eqn, 3rd diagram
\psline[linewidth=.05](33,12.5)(33,7.5)
\psline[linewidth=.02](32.5,12)(36,12)\psline[linewidth=.02](32.5,8)(36,8)
\pscircle*(34,12){.15}\pscircle*(35,12){.15}
\pscircle*(33,11){.15}\pscircle*(33,9){.15}
\rput(34,12.7){\sm{$2^+$}}\rput(35.4,12.7){\sm{$3^-$}}
\rput(32.4,11.3){\sm{$1^+$}}\rput(32.4,9.3){\sm{$1^-$}}
\pscircle*(34,8){.15}\pscircle*(35,8){.15}
\rput(34,7.3){\sm{$2^-$}}\rput(35.4,7.3){\sm{$3^+$}}
\psline[linewidth=.025]{<->}(35.5,8.6)(35.5,11.4)\rput(36,10){$\si$}
\rput(34.2,5.8){$\R D_{\{1\};\{2,3\},\eset}$}
% bottom row
\psline[linewidth=.05](25,3.5)(25,-1.5)\psline[linewidth=.05,arrowsize=6pt]{->}(25,3.5)(25,-.5)
\psline[linewidth=.05](10,3.5)(10,-1.5)\psline[linewidth=.05,arrowsize=6pt]{->}(10,-1.5)(10,2.5)
\rput(25.8,1.5){$3^-$}\rput(9.4,.5){$3^+$}
\psline[linewidth=.05](10,-1.5)(25,-1.5)\psline[linewidth=.05,arrowsize=6pt]{->}(25,3.5)(12,3.5)
\rput(18,-.9){$2^-$}\rput(18,2.9){$2^+$}
\psline[linewidth=.05](10,3.5)(25,3.5)\psline[linewidth=.05,arrowsize=6pt]{->}(10,-1.5)(23,-1.5)
\psline[linewidth=.05](10,3.5)(25,-1.5)\psline[linewidth=.05,arrowsize=6pt]{->}(25,-1.5)(12,2.83)
\psline[linewidth=.05](10,-1.5)(17,.83)
\psline[linewidth=.05](25,3.5)(18,1.17)\psline[linewidth=.05,arrowsize=6pt]{->}(18,1.17)(23,2.83)
\rput(21.2,.5){$1^+$}\rput(14,.5){$1^-$}
\pscircle*(10,3.5){.3}\rput(9.3,3.5){$\eset$}
\pscircle*(10,-1.5){.3}\rput(9.3,-1.5){$3$}
\pscircle*(25,3.5){.3}\rput(25.7,3.5){$2$}
\pscircle*(25,-1.5){.3}\rput(26,-1.5){$23$}
\end{pspicture}
\caption{The first line represents an equivalence of two points in $H^1(\R\ov\cM_{0,2};\Q)$
and a relation between three loops in $H^2(\R\ov\cM_{0,3};\Q)$.
The bottom diagram represents the intersection pattern of the six codimension~2 submanifolds 
$\R D_{I;J,K}\!\approx\!S^1$
with the four real hypersurfaces $\R E_{J',K'}\!\approx\!S^2$ in~$\R\ov\cM_{0,3}$.
The former are labeled by the unique element of $I\!\subset\![3]$
and the sign of $(-1)^{|J|}\!=\!(-1)^{|K|}$;
the latter are labeled by the subset $J',K'\!\subset\![3]$ not containing~1.}
\label{TopolRel_fig}
\end{figure}

\noindent
{\bf Real case with conjugate points only, $\R\ov\cM_{0,\ell}$.}
This orientable compact smooth manifold contains compact 
real hypersurfaces~$\R E_{J,K}$, with $J\!\sqcup\!K\!=\![\ell]$, 
and~$\R H_{J,K}$, with $J\!\sqcup\!K\!=\![\ell]$ and $J,K\!\neq\!\eset$, 
whose generic elements are curves
with two irreducible components that are interchanged by the involution in the first case
and are preserved by~it in the second case.
This manifold also contains compact submanifolds~$\R D_{I;J,K}$,
with \hbox{$I\!\sqcup\!J\!\sqcup\!K\!=\![\ell]$}, $I\!\neq\!\eset$, and $|J|\!+\!|K|\!\ge\!2$, 
of (real) codimension~2 whose generic elements are curves with three components 
two of which are interchanged by the involution;
see Figure~\ref{TopolRel_fig} for some examples. 
The submanifolds~$\R E_{J,K}$ and~$\R D_{I;J,K}$ are orientable, while~$\R H_{J,K}$  
is not if $k\!\ge\!3$.\\

\noindent
By \cite[Theorem~2.2]{RDMhomol},
\begin{enumerate}[label=($\R\C\arabic*$),ref=\arabic*,leftmargin=*]

\item\label{RCgen_it} the submanifolds~$\R E_{J,K}$ and~$\R D_{I;J,K}$ 
generate  $H^*(\R\ov\cM_{0,\ell};\Q)$ as an algebra,

\item\label{RCrel_it} subject only to some expected relations.

\end{enumerate}
The relations arise from
\begin{enumerate}[label=($\R\C\ref{RCrel_it}\alph*$),leftmargin=*]

\item certain submanifolds $\R E_{J,K},\R D_{I;J,K}\!\subset\!\R\ov\cM_{0,\ell}$ 
being disjoint,

\item all real hypersurfaces~$\R E_{J,K}$ in~$\R\ov\cM_{0,2}\!\approx\!\R\P^1$ being points 
and thus equivalent in\linebreak $H^1(\R\ov\cM_{0,2};\Q)$, as indicated by 
the left identity in Figure~\ref{TopolRel_fig}, and 

\item the relation in $H_1(\R\ov\cM_{0,3};\Q)$, or equivalently in $H^2(\R\ov\cM_{0,3};\Q)$,
of Remark~3.4 in~\cite{RealEnum},
as indicated by the right identity and the bottom diagram in Figure~\ref{TopolRel_fig}.

\end{enumerate}
The above relations in $H^1(\R\ov\cM_{0,2};\Q)$ and $H^2(\R\ov\cM_{0,3};\Q)$
pull back to relations in $H^1(\R\ov\cM_{0,\ell};\Q)$ and $H^2(\R\ov\cM_{0,3};\Q)$
via the forgetful morphisms 
$$\R\ov\cM_{0,\ell}\lra \R\ov\cM_{0,2}  \qquad\hbox{and}\qquad
\R\ov\cM_{0,\ell}\lra \R\ov\cM_{0,3}$$
dropping all but two (resp.~three) conjugate pairs of marked points.\\

\noindent
The proof of~($\R\C\ref{RCgen_it}$) in~\cite{RDMhomol} is based
on the presentation of $\R\ov\cM_{0,\ell+1}$ in~\cite{RDMbl} as a sequence 
of blowups of~$\R\ov\cM_{0,\ell}\!\times\!\ov\cM_4$ of three different types
and describes generators for~$H^*(X_{\vr};\Q)$ for every intermediate blowup~$X_{\vr}$, 
assuming~($\R\C\ref{RCgen_it}$) for~$\R\ov\cM_{0,\ell}$.
This approach can also be used to describe generators for~$H_*(\R\ov\cM_{0,\ell};\Z)$
and implies that this module has only 2-torsion.
The proof of~($\R\C\ref{RCrel_it}$) in~\cite{RDMhomol} relates $H^*(\R\ov\cM_{0,\ell+1};\Q)$
to $H^*(\R\ov\cM_{0,\ell'};\Q)$ with $\ell'\!\le\!\ell$ and 
$H^*(\ov\cM_{\ell'};\Q)$ with $\ell'\!\le\!\ell\!+\!1$ by pairing certain cohomology classes 
with the submanifolds~$\R E_{J,K}$ and~$\R D_{I;J,K}$ to show the algebraic independence of 
the former.\\

\begin{figure}
\begin{pspicture}(-2.8,-1.6)(10,1.2)
\psset{unit=.3cm}
% left diagram
\psarc(5,0){5}{90}{270}
\pscircle*(.17,1.29){.2}\pscircle*(1.46,3.54){.2}\pscircle*(3.71,4.83){.2}
\pscircle*(.17,-1.29){.2}\pscircle*(1.46,-3.54){.2}\pscircle*(3.71,-4.83){.2}
\rput(-.5,1.65){\sm{$3^+$}}\rput(.9,3.95){\sm{$2^+$}}\rput(3.5,5.45){\sm{$1^+$}}
\rput(-.6,-1.85){\sm{$3^-$}}\rput(.9,-4.15){\sm{$2^-$}}\rput(3.5,-5.55){\sm{$1^-$}}
\rput(6,0){\begin{Large}$\times$\end{Large}}
\psline(7.5,0)(15.5,0)\pscircle*(8.5,0){.2}\pscircle*(10.5,0){.2}
\pscircle*(12.5,0){.2}\pscircle*(14.5,0){.2}
\rput(8.5,.8){\sm{$4$}}\rput(10.5,-.8){\sm{$3$}}\rput(12.5,.8){\sm{$2$}}\rput(14.5,-.8){\sm{$1$}}
\psline{->}(17.5,0)(22.5,0)
% right diagram
\psarc(30,0){5}{90}{270}
\pscircle*(25.17,1.29){.2}\pscircle*(26.46,3.54){.2}\pscircle*(28.71,4.83){.2}
\pscircle*(25.17,-1.29){.2}\pscircle*(26.46,-3.54){.2}\pscircle*(28.71,-4.83){.2}
\rput(24.5,1.65){\sm{$5^+$}}\rput(28.5,5.45){\sm{$1^+$}}
\rput(24.4,-1.85){\sm{$5^-$}}\rput(28.5,-5.55){\sm{$1^-$}}
\psline(25.46,3.54)(33.46,3.54)\pscircle*(28.46,3.54){.2}
\pscircle*(30.46,3.54){.2}\pscircle*(32.46,3.54){.2}
\rput(28.6,2.74){\sm{$4^+$}}\rput(30.7,4.4){\sm{$3^+$}}\rput(32.7,2.74){\sm{$2^+$}}
\psline(25.46,-3.54)(33.46,-3.54)\pscircle*(28.46,-3.54){.2}
\pscircle*(30.46,-3.54){.2}\pscircle*(32.46,-3.54){.2}
\rput(28.6,-2.74){\sm{$4^-$}}\rput(30.8,-4.4){\sm{$3^-$}}\rput(32.7,-2.74){\sm{$2^-$}}
\psline[linewidth=.025]{<->}(31.46,-2.5)(31.46,2.5)\rput(32.2,0){$\si$}
\psline[linewidth=.025]{<->}(2.46,-2.5)(2.46,2.5)\rput(3.2,0){$\si$}
\end{pspicture}
\caption{The action of the immersion~$\io_{232}^+$ on a typical element in its domain;
each double-headed arrow labeled~$\si$ indicates
the involution on the corresponding real curve.}
\label{ioplus_fig}
\end{figure}
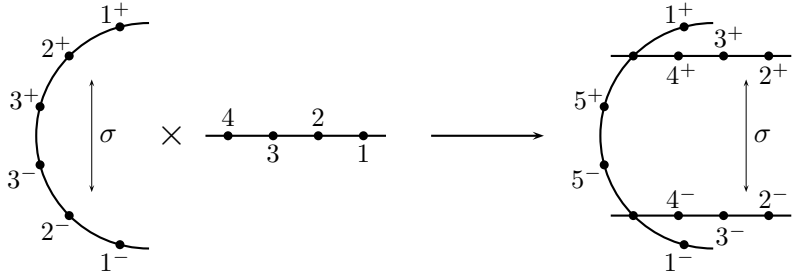

\noindent
Let $R$ be a commutative ring with unity~1 and
$$\cO_{\R\C}(k)=H_*\big(\R\ov\cM_{0,k+1};R\big) \qquad\forall~k\!\in\!\Z^+.$$
For integers $k\!\ge\!1$, $\ell\!\ge\!2$, and $i\!\in\![k]$, let
$$\io_{k\ell i}^+\!: \R\ov\cM_{0,k+1}\!\times\!\ov\cM_{\ell+1}\lra\R\ov\cM_{0,k+\ell}$$
be the immersion identifying the first marked point in the $i$-th conjugate pair of a marked curve 
$[\cC_1]\!\in\!\R\ov\cM_{0,k+1}$ with the last marked point of a marked curve 
$[\cC_2]\!\in\!\ov\cM_{\ell+1}$ and the second marked point in this pair 
with the last marked point of the conjugate marked curve~$[\ov{\cC_2}]$,
keeping the indices of the first $i\!-\!1$ conjugate pairs of marked points of~$[\cC_1]$
the same,
increasing the indices of the last $k\!+\!1\!-\!i$ conjugate pairs of marked points of~$[\cC_1]$
by~$\ell\!-\!1$, and 
turning the first $\ell$~marked points of~$[\cC_2]$ into the first marked of
the conjugate pairs with the indices increased by~$i\!-\!1$;
see Figure~\ref{ioplus_fig} for an example.
This induces homomorphisms
$$\io_{k\ell i*}^+\!: \cO_{\R\C}(k)\!\times\!\cO_{\C}(\ell)\!\approx\!
H_*\big(\R\ov\cM_{0,k+1}\!\times\!\ov\cM_{\ell+1};R\big)\lra
\cO_{\R\C}(k\!+\!\ell\!-\!1).$$
For $k\!\in\!\Z^+$ and $i\!\in\![k]$, define
$$\io_{k1i*}^+\!: \cO_{\R\C}(k)\!\times\!\cO_{\C}(1)\lra\cO_{\R\C}(k), 
\qquad \io_{k1i*}^+(\eta,r)=r\eta.$$
For $k,\ell_1,\ldots,\ell_k\!\in\!\Z^+$, define
\begin{gather*}
\circ\!:\cO_{\R\C}(k)\!\times\!\cO_{\C}(\ell_1)\!\times\!\ldots\!\times\!\cO_{\C}(\ell_k)\lra
\cO_{\R\C}\big(\ell_1\!+\!\ldots\!+\!\ell_k\big),\\
\eta\!\circ\!(\eta_1,\ldots,\eta_k)
=\io_{(1+\ell_2+\ldots+\ell_k)\ell_11*}^+
\big(\ldots\io_{(k-1+\ell_k)\ell_{k-1}(k-1)*}^+\big(\io_{k\ell_kk*}^+(\eta,\eta_k),\eta_{k-1}\big),
\ldots,\eta_1\big). 
\end{gather*}
This determines a right module structure on $(\cO_{\R\C}(k))_{k\in\Z^+}$ over
the operad $(\!(\cO_{\C}(\ell)\!)_{\ell\in\Z^+},\circ)$.

\section{Open Problems}

\noindent
{\bf Real case with conjugate points only, $\R\ov\cM_{0,\ell}$.}
The approach of~\cite{RDMbl,RDMhomol} to~($\R\C\ref{RCgen_it}$) on page~\pageref{RCgen_it}
implies that $H^*(\R\ov\cM_{0,\ell};\Z)$ is linearly generated by orientable intersections
of the submanifolds $\R E_{J,K}$, $\R H_{J,K}$, and~$\R D_{I;J,K}$
and contains only 2-torsion.
For example, the torsion of $H_1(\R\ov\cM_{0,3};\Z)$ is generated by 
$\R H_{\{1,2\},\{3\}}\!\cap\!\R H_{\{1\},\{2,3\}}$, as shown in \cite[Section~3]{RealEnum}.
However, it still remains to determine the amount of 2-torsion and 
the ring structure of~$H^*(\R\ov\cM_{0,\ell};\Z)$ for $\ell\!\ge\!4$.
The algebraic generators for the torsion of~$H_*(\R\ov\cM_{0,\ell};\Z)$ should correspond to trees,
with the univalent and bivalent vertices labeled by nonempty subsets of~$[\ell]$ forming
a partition of~$[\ell]$.
A more algebraic question is to describe the structure of
the right module structure~$(\cO_{\R\C}(k))_{k\in\Z^+}$ over
the operad $(\!(\cO_{\C}(\ell)\!)_{\ell\in\Z^+},\circ)$ explicitly,
in the spirit of~\cite{Getz95} and \cite[Theorem~2.14]{EHKR}.\\

\noindent
{\bf The general real case, $\R\ov\cM_{k,\ell}$ with $k,\ell\!\in\!\Z^+$.}
The approach of~\cite{RDMbl,RDMhomol} to~($\R\C\ref{RCgen_it}$) also implies that 
$H_*(\R\ov\cM_{k,\ell};\Z)$ with $k\!\ge\!1$ is linearly generated by orientable intersections
of the submanifolds $\R H_{J,K}$ and~$\R D_{I;J,K}$ and contains only 2-torsion.
The former are now indexed by pairs $J\!\equiv\!(J_{\R},J_{\C})$ and $K\!\equiv\!(K_{\R},K_{\C})$~with 
$$[k]=J_{\R}\!\sqcup\!K_{\R}, \qquad [\ell]=J_{\C}\!\sqcup\!K_{\C}, 
\qquad\hbox{and}\quad 
\big|J_{\R}\big|\!+\!2\big|J_{\C}\big|,\big|K_{\R}\big|\!+\!2\big|K_{\C}\big|\ge2\,,$$
and the latter by the partitions $[\ell]\!=\!I\!\sqcup\!J\!\sqcup\!K$ with $|J|\!+\!|K|\!\ge\!2$.
Similarly, $H^*(\R\ov\cM_{k,\ell};\Z)$ is generated by some linear combinations of these intersections.
In fact, the most delicate type of blowups appearing in~\cite{RDMbl,RDMhomol} no longer enters
into consideration.
However, it remains to determine the amount of 2-torsion in~$H_*(\R\ov\cM_{k,\ell};\Z)$
and the ring structure even of~$H^*(\R\ov\cM_{k,\ell};\Q)$.\\

\noindent
For $k\!\in\!\Z^+$ and $\ell\!\in\!\Z^{\ge0}$, let
$$\cO_{\R}(k,\ell)=\begin{cases}\Q,&\hbox{if}~k\!+\!\ell\!=\!1;\\
H_*(\R\ov\cM_{k+1,\ell};\Q),&\hbox{if}~k\!+\!\ell\!\ge\!2;\end{cases}
\qquad \cO_{\R}(k)=\bigoplus_{\ell=0}^{\i} \cO_{\R}(k,\ell).$$
Similarly to the previous cases, there are node-identifying immersions
\begin{equation*}\begin{split}
\io_{k\ell k'\ell' i}^{\R}\!: \R\ov\cM_{k+1,\ell}\!\times\!\R\ov\cM_{k'+1,\ell'}
&\lra \R\ov\cM_{k+k',\ell+\ell'} \qquad\hbox{and}\\
\io_{k\ell\ell' i}^+\!: \R\ov\cM_{k+1,\ell}\!\times\!\ov\cM_{\ell'+1}
&\lra \R\ov\cM_{k+1,\ell+\ell'-1},
\end{split}\end{equation*}
with $k\!+\!\ell\!\ge\!2$ in both cases, 
$i\!\in\![k\!+\!1]$ and $k'\!+\!\ell'\!\ge\!2$ in the first case, 
and $i\!\in\![\ell]$ and $\ell'\!\ge\!2$ in the last case.
They induce homomorphisms
\begin{equation*}\begin{split}
\io_{k\ell k'\ell' i*}^{\R}\!: \cO_{\R}(k,\ell)\!\times\!\cO_{\R}(k',\ell')\!\approx\!
H_*\big(\R\ov\cM_{k+1,\ell}\!\times\!\R\ov\cM_{k'+1,\ell'};\Q\big)
&\lra \cO_{\R}\big(k\!+\!k'\!-\!1,\ell\!+\!\ell'\big) \quad\hbox{and}\\
\io_{k\ell\ell' i*}^+\!: \cO_{\R}(k,\ell)\!\times\!\cO_{\C}(\ell')\!\approx\!
H_*\big(\R\ov\cM_{k+1,\ell}\!\times\!\ov\cM_{\ell'+1};\Q\big)
&\lra\cO_{\R}\big(k,\ell\!+\!\ell'\!-\!1\big).
\end{split}\end{equation*}
We define the homomorphisms
\begin{gather*}
\io_{k\ell10i*}^{\R}\!: \cO_{\R}(k,\ell)\!\times\!\cO_{\R}(1,0)\lra\cO_{\C}(k,\ell), \qquad
\io_{10k'\ell' 1*}^{\R}\!: \cO_{\R}(1,0)\!\times\!\cO_{\R}(k',\ell')\lra\cO_{\C}(k',\ell'),\\
\hbox{and}\qquad
\io_{k\ell 1 i*}^+\!: \cO_{\R}(k,\ell)\!\times\!\cO_{\C}(1)\lra\cO_{\C}(k,\ell)
\end{gather*}
similarly to~\eref{prodcases_e}.
These bilinear maps induce multilinear maps
\begin{equation*}\begin{split}
\circ\!:\cO_{\R}(k,\ell)\!\times\!
\cO_{\R}(m_1,\ell_1')\!\times\!\ldots\!\times\!\cO_{\R}(m_k,\ell_k')&\lra
\cO_{\R}\big(m_1\!+\!\ldots\!+\!m_k,\ell\!+\!\ell_1'\!+\!\ldots\!+\!\ell_k'\big) \qquad\hbox{and}\\
\circ\!:\cO_{\R}(k,\ell)\!\times\!\cO_{\C}(m_1)\!\times\!\ldots\!\times\!\cO_{\C}(m_{\ell})&\lra
\cO_{\R}\big(k,m_1\!+\!\ldots\!+\!m_{\ell}\big)
\end{split}\end{equation*}
similarly to~\eref{expandcomp_e}.
The last two maps determine an operad $(\cO_{\R}(k))_{k\in\Z^+}$
of graded vector spaces and a right module structure on $(\cO_{\R}(k,\ell)\!)_{\ell\in\Z^+}$
over the operad $(\!(\cO_{\C}(\ell)\!)_{\ell\in\Z^+},\circ)$ for each $k\!\in\!\Z^+$.\\

\noindent
As suggested by A.~Voronov, the above structures can be combined into an object resembling 
a bicolored operad; see \cite[Section~1.1]{Gir}. 
For $n\!\in\!\Z^+$, let
$$\cO_{\R}'(n)=\cO_{\C}(n)\sqcup\bigsqcup_{\begin{subarray}{c}k+\ell=n\\ 
k\ge0,\ell\ge1 \end{subarray}}\!\!\!\!\!\cO_{\R}(k\!+\!1,\ell\!-\!1)\,.$$
For each fixed $n\!\in\!\Z^+$, auxiliary grading is provided by $k\!\in\!\Z^{\ge0}$ above,
with the elements of~$\cO_{\C}(n)$ being of the auxiliary degree~0.
For $x\!\in\!\cO_{\R}'(n)$, let $|x|\!=\!n$ and
$$|x|_+=\begin{cases}n,&\hbox{if}~x\!\in\!\cO_{\C}(n);\\
\ell\!-\!1,&\hbox{if}~x\!\in\!\cO_{\R}(k\!+\!1,\ell\!-\!1).\end{cases}$$
We take the set of colors to be $\cC\!\equiv\!\{+,\R\}$ and define
$$\tn{out}\!:\cO_{\R}'(n)\lra\cC, \qquad 
\tn{out}(x)=\begin{cases}+,&\hbox{if}~x\!\in\!\cO_{\C}(n);\\
\R,&\hbox{if}~x\!\in\!\cO_{\R}(k\!+\!1,\ell\!-\!1).\end{cases}$$
For $i\!\in\![n]$, define
$$\tn{in}_i\!:\cO_{\R}'(n)\lra\cC, \qquad 
\tn{in}_i(x)=\begin{cases}+,&\hbox{if}~i\!\le\!|x|_+;\\
\R,&\hbox{if}~i\!>\!|x|_+.
\end{cases}$$
Let $1_+\!=\!1\!\in\!\cO_{\C}(1)$ and $1_{\R}\!=\!1\!\in\!\cO_{\R}(1,0)$.
The above auxiliary grading is preserved by the~maps
$$\circ_i\!:\big\{(x,y)\!\in\!\cO_{\R}'(m)\!\times\!\cO_{\R}'(n)\!:
\tn{in}_i(x)\!=\!\tn{out}(y)\!\big\}\lra \cO_{\R}'(m\!+\!n\!-\!i), 
\quad i\!\in\![m],$$
induced by the homomorphisms~$\io_{k\ell k'\ell' i*}^{\R}$ and $\io_{k\ell\ell' i*}^+$
defined in the previous paragraph.
These maps satisfy a modification of~\cite[(1.1.4a)]{Gir},
\BE{114a_e} x\!\circ_i\!(y\!\circ_j\!z)=\begin{cases}
(x\!\circ_i\!y)\!\circ_{j+i-1}\!z,&\hbox{if}~\tn{out}(y)\!=\!\tn{out}(z);\\
(x\!\circ_i\!y)\!\circ_{j+|x|_+}\!z,&\hbox{if}~\tn{out}(y)\!\neq\!\tn{out}(z);
\end{cases}  \quad\hbox{if}~1\le i\le |x|,\,1\le j\le |y|.\EE
These maps also satisfy a modification of~\cite[(1.1.4b)]{Gir},
\BE{114b_e} 
(x\!\circ_i\!y)\!\circ_{j+|y|-1}\!z=\begin{cases}
(x\!\circ_j\!z)\!\circ_i\!y,&\hbox{if}~\tn{out}(y)\!=\!+;\\
(x\!\circ_j\!z)\!\circ_{i+|z|_+}\!y,&\hbox{if}~\tn{out}(y)\!=\!\R;
\end{cases}  \quad\hbox{if}~1\le i< j\le |x|.\EE
The identities~\eref{114a_e} and~\eref{114b_e} hold whenever both sides are defined.
They agree with \cite[(1.1.4a)]{Gir} and \cite[(1.1.4b)]{Gir}
whenever~$\tn{out}(y)\!=\!\tn{out}(x)$ and $\tn{out}(y)\!=\!\R$, respectively.
The identity \cite[(1.1.4c)]{Gir} holds without any modifications.\\

\noindent
A natural question is to describe the above algebraic structures explicitly
and perhaps to find a more succinct way of presenting them.\\

\vspace{.2in}

\noindent
{\bf Acknowledgment.} The author would like to thank A.~Voronov for discussions
on operads.\\

\vspace{.2in}

\noindent
{\it Department of Mathematics, Stony Brook University, Stony Brook, NY 11794\\
azinger@math.stonybrook.edu}\\


\begin{thebibliography}{99}

\bibitem{RDMhomol} X.~Chen, P.~Georgieva, and A.~Zinger,
{\it The cohomology ring of the Deligne-Mumford space
of real rational curves with conjugate marked points}, math/2305.08798

\bibitem{RDMbl} X.~Chen and A.~Zinger, {\it Blowdowns of the Deligne-Mumford spaces
of real rational curves}, math/2305.08811

\bibitem{EHKR} P.~Etingof, A.~Henriques, J.~Kamnitzer, and E.~Rains,
{\it The cohomology ring of the real locus of the moduli space of stable curves 
of genus~0 with marked points}, Ann.~of Math.~171 (2010), no.~2, 731-–777

\bibitem{RealEnum} P.~Georgieva and A.~Zinger,
{\it Enumeration of real curves in $\C\P^{2n-1}$ 
and a WDVV relation for real Gromov-Witten invariants}, 
Duke Math.~J.~166 (2017), no.~17, 3291-–3347

\bibitem{Getz95} E.~Getzler, {\it Operads and moduli spaces of genus 0 Riemann surfaces}, 
in {\it The Moduli Space of Curves}, Progr.~Math.~129 (1995), Birkh\"{a}user,  199–-230

\bibitem{Gir} S.~Giraudo,
{\it Colored operads, series on colored operads, and combinatorial generating systems},
Discrete Math.~342 (2019), no.~6, 1624–-1657

\bibitem{Keel} S.~Keel, 
{\it Intersection theory of moduli spaces of stable $n$-pointed curves of genus zero},
Trans.~AMS 330 (1992), no.~2, 545--574 

\bibitem{Krasnov} V.~Krasnov, {\it Real algebraic maximal varieties}, 
Mat.~Zametki 73 (2003), 853–-860

\bibitem{Rains} E.~Rains, {\it The homology of real subspace arrangements},
J.~Topol.~3 (2010), no.~4, 786-–818


\end{thebibliography}
\end{document}